\documentstyle{amsart}
\newcommand{\genus}{\operatorname{genus}}
\newcommand{\Stab}{\operatorname{Stab}}
\newcommand{\Aut}{\operatorname{Aut}}
\newcommand{\Var}{\operatorname{Var}}
\newtheorem{theorem}{Theorem}
\newtheorem{corollary}{Corollary}
\newtheorem{lemma}{Lemma}
\begin{document}
\title[Commutators and Squares]
        {Products of Commutators and Products of Squares in a Free Group}
\author[L. P. Comerford]{Leo P. Comerford, Jr.}
\address{Department of Mathematics\\Eastern Illinois University\\
        Charleston, Illinois 61920\\U. S. A.}
\email{cflpc@@eiu.edu}
\author[C. C. Edmunds]{Charles C. Edmunds}
\address{Department of Mathematics\\Mount Saint Vincent University\\
        Halifax, Nova Scotia B3M 2J6\\Canada}
\email{CEDMUNDS@@linden.msvu.ca}
\maketitle
\begin{abstract}
A classification of the ways in which an element of a free group can be 
expressed as a product of commutators or as a product of squares is given.  
This is then applied to some particular classes of elements.  Finally, a 
question about expressing a commutator as a product of squares is addressed.
\end{abstract}

\section{Introduction}\label{se:intro}

Our first aim is to provide a characterization of the set of all solutions to 
the equations 
\begin{equation}
[x_1,y_1]\dots [x_g,y_g]=U\label{eq:or}
\end{equation}
\begin{equation}
x_1^2\dots x_g^2=U\label{eq:nonor}
\end{equation}
in a free group, where $U$ is a product of $g$ but not fewer commutators in 
\eqref{eq:or} and $U$ is a product of $g$ but not fewer squares and is not 
a product of fewer 
than $g/2$ commutators in \eqref{eq:nonor}.  Our description of the set of 
solutions, which we give in \S\ref{se:theorems}, is based on work of 
C.~Edmunds \cite{E-Endo2}. It provides an algebraic version of a result 
obtained by M.~Culler \cite[Theorem 4.1]{Cul81} by topological methods. This 
characterization of solutions is somewhat more explicit than those given by 
L.~Comerford and Edmunds \cite{CE-Sing} and by R.~Grigorchuk and P.~Kurchanov 
\cite{GrigKurch92} for a larger class of equations, and seems easier to use to 
get descriptions of solutions to particular equations.  We use our 
method to classify solutions to certain classes of equations in 
\S\ref{se:applics}.  In \S\ref{se:comms} we give another application, this 
time to expressions of commutators as products of squares in a free group.

\section{Products of Commutators, Products of Squares}\label{se:theorems}

We begin by fixing some notation and terminology.  Let $H$ be the free group 
on $A=\{a_1,\dots,b_1,\dots,c_1,\dots\}$ and let $F$ be the free group on 
$X=\{x_1,\dots,y_1,\dots\}$.  We call elements of $A$ {\em constants}, 
elements of $X$ {\em variables}, and elements of $X\cup X^{-1}$ {\em letters}.  
Length of elements of $H$ or $F$ relative to these generating sets is denoted 
by $|\cdot |$.  For $W\in F$, we let $\Var(W)$ be the set of variables 
occurring in $W$ and let $|\Var(W)|$ be its cardinality.

We denote by $G'$ the commutator subgroup of a group $G$, and by $2G$ the 
subgroup generated by squares of elements of $G$.  For $U\in G'$, we let 
$\genus^+(U)$ be the minimal number of commutators of which $U$ is a product, 
and for $U\in 2G$, we let $\genus^-(U)$ be the minimal number of squares of 
which $U$ is a product.  We set $\genus^+(1)=\genus^-(1)=0$, and let 
$\genus^+(U)=\infty$ if $U\not\in G'$ and let $\genus^-(U)=\infty$ if 
$U\not\in 2G$.

Since $G'\subseteq 2G$ and since there are automorphisms of the free group on 
$x, y, z$ sending $x^2[y,z]$ and $x^2y^2z^2$ to one another, it follows that 
for $U\in G'$, $\genus^-(U)\le 2\genus^+(U)+1$.  There is, however, 
no general upper bound on $\genus^+(U)$ in terms of $\genus^-(U)$ for 
$U\in G'$.  For instance, if $U=([a_1,b_1]\dots[a_g,b_g])^2$, $\genus^-(U)=1$ 
but $\genus^+(U)=2g$ (cf. \cite{CE-Genus}).

An element $W$ of $F$ is called {\em quadratic} if each variable that occurs 
in $W$ occurs exactly twice, with exponents $+1$ or $-1$.  We call a quadratic 
word $W$ {\em orientable} if $W\in F'$, and {\em nonorientable} otherwise.  
Note that every quadratic word is an element of $2F$.  A quadratic word $W$ is 
{\em irredundant} if there is no pair of distinct, noninverse letters 
$x,y$ which appear in $W$ only in subwords $(xy)^{\pm 1}$.                    

With any quadratic word $W$ we may associate a closed surface $S_W$ by writing 
$W$ around the boundary of a disk and identifying edges labeled by the same 
variable, respecting orientation.  Note that the orientability of $S_W$ is the 
same as that of $W$.   We define the {\em Euler 
characteristic} of $W$, $\chi(W)$, to be that of $S_W$.  By classical results, 
$\chi(W)=2-2\genus^+(W)$ if $W$ is orientable and $\chi(W)=2-\genus^-(W)$ if 
$W$ is nonorientable.  Also, if $\alpha$ is an automorphism of $F$ and both 
$W$ and $W\alpha$ are quadratic, $\chi(W)=\chi(W\alpha)$.  If $W$ is an 
orientable quadratic word with $\genus^+(W)=g$, there is an automorphism 
$\alpha$ of $F$ with $W\alpha=[x_1,y_1]\dots[x_g,y_g]$, and if $W$ is a 
nonorientable quadratic word with $\genus^-(W)=g$, there is an automorphism 
$\alpha$ of $F$ with $W\alpha=x_1^2\dots x_g^2$.

A {\em solution} to an equation 
\begin{equation}
W(x_1,\dots)=U(a_1,\dots)\label{eq:general}
\end{equation}
is a homomorphism $\phi:F\to H$ such that $W\phi=U$.  A solution $\phi$ to 
\eqref{eq:general} is called {\em cancellation-free} of $x\phi\ne 1$ for each 
variable $x$ occuring in $W$ and if $W(x_1\phi,\dots)\equiv U(a_1,\dots)$, 
where $\equiv$ denotes equality in the free semigroup on $A\cup A^{-1}$.

Finally, we need to say something about stabilizers of elements of $F$.  For 
$W\in F$, we define $F_W$ to be the subgroup of $F$ generated by $\Var(W)$, 
we let $\Stab_F(W)=\{\alpha\in\Aut(F) : W\alpha=W\}$,   
and we let $\Stab_{F_W}=\{\alpha\in\Aut(F_W) : W\alpha=W\}$.  From a result of 
G.~Rosenberger \cite[Theorem 2.2]{Ros89}, it follows that if 
$W=[x_1,y_1]\dots[x_g,y_g]$ or $W=x_1^2\dots x_g^2$ and if 
$\alpha\in\Stab_F(W)$, then $\alpha$ maps $F_W$ onto $F_W$.  Thus,
for these values of $W$, if 
$\alpha\in\Aut(F)$ and $W\alpha=W$, the restriction of $\alpha$ to $F_W$ is an 
element of $\Aut(F_W)$, so we need not distinguish between the stabilizers of 
$W$ in $F$ and in $F_W$, and we denote both by $\Stab(W)$.

Note that if $\phi$ is a solution to \eqref{eq:general} and 
$\sigma\in\Stab_F(W)$, then $\sigma\phi$ is also a solution to 
\eqref{eq:general}.  We say that two solutions $\phi_1$ and $\phi_2$ to 
\eqref{eq:general} are in the same {\em stabilizer class} if 
$\phi_2=\sigma\phi_1$ for some $\sigma\in\Stab_F(W)$.

As general references, and for unexplained notation or terminology, we refer 
the reader to the books of R.~Lyndon and P.~Schupp \cite{LS} and of W.~Massey 
\cite{Mas}.

We are now ready to state our classification results, for the orientable and 
the nonorientable cases.

\begin{theorem}\label{th:main}
Let $U$ be a nontrivial element of the free group $H$.
\begin{itemize}
\item[(a)] 
If $\genus^+(U)=g<\infty$, and $\phi$ is a solution to
\begin{equation}
[x_1,y_1]\dots[x_g,y_g]=U,\tag{\ref{eq:or}}
\end{equation}
there is an irredundant orientable quadratic word $W\in F$ with 
$\genus^+(W)=g$ and a cancellation-free solution $\psi$ to $W=U$ such that for 
any automorphism $\gamma_W$ of $F$ with $W\gamma_W=[x_1,y_1]\dots[x_g,y_g]$, 
there is an $\alpha\in\Stab([x_1,y_1]\dots[x_g,y_g])$ such that 
$\phi=\alpha\gamma_W^{-1}\psi$.
\item[(b)]
If $\genus^-(U)=g<\infty$, $\genus^+(U)\ge g/2$, and $\phi$ is a solution to
\begin{equation}
x_1^2\dots x_g^2=U,\tag{\ref{eq:nonor}}
\end{equation}
then there is an irredundant nonorientable quadratic word $W\in F$ with 
$\genus^-(W)=g$ and a cancellation-free solution $\psi$ to $W=U$ such that for 
any automorphism $\gamma_W$ of $F$ with $W\gamma_W=x_1^2\dots x_g^2$, there is 
an $\alpha\in\Stab(x_1^2\dots x_g^2)$ such that 
$\phi=\alpha\gamma_W^{-1}\psi$.
\end{itemize}
\end{theorem}

We note that in this theorem we may assume, by composing with an inner 
automorphism if necessary, that $U$ is cyclically reduced.

In viewing $U$ as an image of a quadratic word $W$, it is convenient to take 
$U$ and $W$ to be cyclic words, that is, words written around a circle.  If 
$U$ is an image of $W$ under a map $\psi$ as ordinary words, this remains the 
case if we view $U$ and $W$ as cyclic words.  Our goal is to express $U$ as a 
cancellation-free image of a quadratic word $W$ under a map $\psi$.  If this 
is the case with $U$ and $W$ cyclic words, we may be required to split a 
variable in $W$ to obtain such a representation as ordinary words.  For 
example, $U=a^{-1}c^{-1}abcb^{-1}$ is a cancellation-free image of 
$W=x^{-1}y^{-1}xy$ under $\psi : x\mapsto ab, y\mapsto c$ as cyclic words, but 
$U$ is not a cancellation-free image of $W$ as ordinary words.  Since the 
ordinary word $U$ begins within the image of $x^{-1}$, we replace $x$ by 
$x_1x_2$ and define $(x_1)\psi=a$ and $(x_2)\psi=b$.  We now find that as 
ordinary words, $U$ is a cancellation-free image under $\psi$ of 
$W'=x_1^{-1}y^{-1}x_1x_2yx_2^{-1}$.  Notice that in this process if $W$ is 
irredundant and orientable or nonorientable as a cyclic word, the same will be 
true of the ordinary word $W'$.

For a quadratic word $W$, the set of all irredundant quadratic cyclic words 
with the genus and orientability of $W$, distinct up to automorphisms of $F$ 
that permute $X\cup X^{-1}$, is called the set of {\em Wicks forms} for $W$.  
M.~Wicks showed in \cite{Wicks62} that the orientable Wicks forms of genus one 
are $x^{-1}y^{-1}xy$ and $x^{-1}y^{-1}z^{-1}xyz$, and in \cite{Wicks73} that 
the nonorientable Wicks forms of genus two are $x^2y^2$, $xy^{-1}xy$, 
$z^{-1}x^2zy^2$, and $xzxyz^{-1}y$.  The nine maximal length orientable Wicks 
forms of genus two are listed by J.~Comerford, L.~Comerford, and Edmunds in 
\cite{CCE}.  A.~Vdovina [written communication] has produced lists of the 
maximal length Wicks forms for nonorientable genus three and four and for 
orientable genus three. 

Edmunds showed in \cite{E-Endo2} that if $Q\in F$, $Q$ is quadratic, $U\in H$, 
and $\phi$ is a solution to $Q=U$, then there is an endomorphism $\beta$ of 
$F$ such that $Q\beta$ is an irredundant quadratic word and $Q\beta=U$ has a 
cancellation-free solution.  We modify Edmunds' proof to show that with 
$Q=[x_1,y_1]\dots[x_g,y_g]$ or $Q=x_1^2\dots x_g^2$ and our hypotheses on $U$, 
we may arrange to have $\beta$ an automorphism of $F$ and 
$\psi=\beta^{-1}\phi$.  This will establish the conclusions of the theorem.

\begin{pf*}{Proof of Theorem \ref{th:main}}

We initially set $W=[x_1,y_1]\dots[x_g,y_g]$ or $W=x_1^2\dots x_g^2$ and
$\psi=\phi$, and proceed by induction first on 
$N(W,\psi)=\sum_{x\in\Var(W)}|x\psi|$ and second on $|\Var(W)|$.  We shall 
show that if $W$, $\psi$ do not satisfy the conclusions of the theorem, there 
is an automorphism $\beta$ of $F$ such that $W\beta$ is quadratic and either 
$N(W\beta,\beta^{-1}\psi)<N(W,\psi)$ or $N(W\beta,\beta^{-1}\psi)=N(W,\psi)$ 
and $|\Var(W\beta)|<|\Var(W)|$; we then replace $W$ by $W\beta$ and $\psi$ by 
$\beta^{-1}\psi$ and appeal to the induction hypothesis.  Note that if $W$, 
$\psi$ do not satisfy the conclusions of the theorem, then either $W$ is 
redundant (which is not the case initially, but could occur later in the 
process), there are letters $x$, $y$ such that $xy$ is a subword of $W$ and 
there is cancellation in the product $(x\psi)(y\psi)$, or there is a variable 
$x$ in $W$ such that $x\psi=1$.
 
First, if $W$ is redundant, 
there are distinct letters $x$ and $y$ that occur in $W$ only in subwords
$(xy)^{\pm 1}$.  We find that if $\beta : x \mapsto xy^{-1}$, 
$N(W\beta,\beta^{-1}\psi)\le N(W,\psi)$ and $|\Var(W\beta)|<|\Var(W)|$. 
(When defining an endomorphism of $F$, we take it to fix all letters whose 
images are not specified.)

Suppose that there are letters $x$ and $y$ such that $xy$ is a subword of $W$ 
and there is cancellation in the product $(x\psi)(y\psi)$.  Let 
$x\psi\equiv AB$ and $y\psi\equiv B^{-1}C$ with $B\ne 1$ and the product $AC$ 
freely reduced.  In this 
case we define $\beta$ by $x\beta=xz$, $y\beta=z^{-1}y$ where $z$ is a 
variable that has not previously appeared in any quadratic word used in our 
process, and we specify that $z\psi=B$.  There is no harm in this last 
requirement, since images under $\psi$ of variables not in $W$ are irrelevant 
to the value of $W\psi$.  Here we see that 
$N(W\beta,\beta^{-1}\psi)<N(W,\psi)$.

Finally, suppose that $x\psi=1$ for some variable $x$ occurring in $W$.  We 
must show that there is an automorphism $\beta$ of $F$ with $W\beta=W\tau$, 
where $\tau$ is the endomorphism of $F$ defined by $x\tau=1$.  Note that our 
hypotheses ensure that, since $W\tau=U$ has a solution, 
$\genus^+(W\tau)=\genus^+(W)$ if $W$ is orientable, and that if $W$ is 
nonorientable, then $W\tau$ is nonorientable as well and 
$\genus^-(W\tau)=\genus^-(W)$.  (We use the fact that if $W$ is a 
nonorientable quadratic word with $\genus^-(W)=g$ and if $V$ is an orientable 
quadratic word which is a homomorphic image of $W$, then 
$\genus^+(V)\le (g-1)/2$; cf.~\cite[Proposition I.6.10]{LS}.)  
It follows that the initial and terminal vertices 
of the edge $e_x$ labeled by $x$ on $S_W$ are distinct.  Let $v$ be the 
initial 
vertex of this edge.  We define $\beta$ so that if $y$ is a variable other 
than $x$ with $e_y$ having initial but not terminal vertex at $v$ then 
$y\beta=xy$, so that if $y$ is a variable with $e_y$ having terminal but not 
initial vertex at $v$ then $y\beta=yx^{-1}$, and so that if $y$ is a variable 
with $e_y$ having both initial and terminal vertices at $v$ then 
$y\beta=xyx^{-1}$.  Now $\beta$ is an automorphism, in fact a Whitehead 
automorphism, of $F$ and one may check that $W\beta=W\tau$.  We see, then, 
that $N(W\beta,\beta^{-1}\psi)\le N(W,\psi)$ and $|\Var(W\beta)|<|\Var(W)|$.
\end{pf*}

Theorem \ref{th:main} shows us, then, that we can get 
representatives of all stabilizer classes of solutions to \eqref{eq:or} and 
\eqref{eq:nonor} by finding all the ways in which the right-hand side is a 
cancellation-free image of a Wicks form for the left-hand side.  Note that the 
choice of $\gamma_W$ does not affect the stabilizer class, for if $\psi$ is a 
cancellation-free solution to $W=U$ with $W$ a Wicks form for 
$Q=[x_1,y_1]\dots[x_g,y_g]$ or $Q=x_1^2\dots x_g^2$, and if $\gamma_1$ and 
$\gamma_2$ are two automorphisms of $F$ sending $W$ to $Q$, then 
$\gamma_1^{-1}\gamma_2\in\Stab(Q)$ and so $\gamma_1^{-1}\psi$ and 
$\gamma_2^{-1}\psi$ are in the same stabilizer class of solutions to 
\eqref{eq:or} or \eqref{eq:nonor}.  Also note that if $W_1$ and $W_2$ are 
Wicks forms for $Q$ and $\psi_1$ and $\psi_2$ are cancellation-free solutions 
to $W_1=U$ and $W_2=U$ respectively and if there is a relabelling of variables 
$\rho$ (that is, an automorphism $\rho$ of $F$ that permutes $X\cup X^{-1}$) 
such that $W_2=W_1\rho$ and $\psi_2=\rho^{-1}\psi_1$, then for any maps 
$\gamma_{W_1}$ and $\gamma_{W_2}$, $\gamma_{W_1}^{-1}\psi_1$ and 
$\gamma_{W_2}^{-1}\psi_2$ are in the same stabilizer class of solutions to 
\eqref{eq:or} or \eqref{eq:nonor}.  Thus we find representatives of all 
stabilizer classes of solutions to \eqref{eq:or} or \eqref{eq:nonor} by 
finding all the ways the right-hand side is a cancellation-free image of a 
Wicks form for the left-hand side, up to relabeling of variables in the Wicks 
form.

Now suppose that $W\in F$ is a Wicks form with $|W|=k$ and that $U\in H$ is 
cyclically reduced and $|U|=n$.  The number of ways in which $U$ might be a 
cancellation-free image of $W$, as cyclic words, is bounded by the number of 
ways in which $U$ can be factored as a product of $k$ subwords, 
${{n+k}\choose{k}}$, times the number of ways to make a first identification 
of a letter of $W$ with a chosen subword of $U$, which is $k$.  Now 
$k{{n+k}\choose{k}}$ is a polynomial of degree $k$ in $n$ and, since every 
vertex on $S_W$ has degree at least three, we find that $k\le 6(1-\chi(W))$.  
Further, the number of orientable or nonorientable Wicks forms of a given 
genus is finite.  Thus we have the following, which is an instance of a more 
general result of Grigorchuck and Kurchanov \cite[Theorem 4]{GrigKurch92}.

\begin{corollary}
Let $Q=[x_1,y_1]\dots[x_g,y_g]$ or $Q=x_1^2\dots x_g^2$ and let $U\in H$.  
Suppose that $U\in H'$ and $\genus^+(U)=g$ if $Q=[x_1,y_1]\dots[x_g,y_g]$, and 
that $U\in 2H$ with $\genus^-(U)=g$ and 
$\genus^+(U)\ge g/2$ if $Q=x_1^2\dots x_g^2$.  There is an algorithm to 
compute a complete set of representatives for stabilizer classes of solutions 
to $Q=U$ whose number of steps is bounded by a polynomial in $|U|$ of degree 
$12g-6$ if $Q=[x_1,y_1]\dots[x_g,y_g]$ and of degree $6g-6$ if 
$Q=x_1^2\dots x_g^2$ and $g>1$.
\end{corollary}

\section{Applications}\label{se:applics}

We use our classification of solutions to sharpen a theorem of R.~Burns, 
Edmunds, and I.~Farouqi 
\cite[Theorem 1]{BEF}, which in turn improved upon a result of 
Ju.~Hmelevski\v{\i} \cite{Hmel71}.  All of this is based on pioneering work of 
A.~I.~Mal'cev \cite{Mal62}.

\begin{theorem}\label{th:BEF}
Suppose that $U\in H$ is nontrivial and cyclically reduced and that there are 
solutions to
\begin{equation}
[x,y]=U.\label{eq:comm}
\end{equation}
Then every solution to \eqref{eq:comm} is in the stabilizer class of a 
solution $\phi_0$ with $|x\phi_0|\le(1/2)|U|$, $|y\phi_0|\le(1/2)|U|$, and 
$|x\phi_0|+|y\phi_0|\le|U|-1$.  Further, $U$ has a cyclic permutation $U^*$ 
such that every solution to $[x,y]=U^*$ is in the stabilizer class of a 
solution $\phi_0^*$ with $|x\phi_0^*|\le(1/2)|U|-1$, 
$|y\phi_0^*|\le(1/2)|U|-1$, and $|x\phi_0^*|+|y\phi_0^*|\le(2/3)|U|$.
\end{theorem}

\begin{pf*}{Proof of Theorem \ref{th:BEF}}

Given the orientable Wicks forms of genus one, 
we know that as a cyclic word $U$ factors as either
$X^{-1}Y^{-1}XY$ or $X^{-1}Y^{-1}Z^{-1}XYZ$ with $X$, $Y$, $Z$ nonempty.  As an 
ordinary word, then, we find that up to change of notation either
\begin{equation}
U\equiv X_1^{-1}Y^{-1}X_1 X_2 Y X_2^{-1}\label{eq:twovar}
\end{equation}
or
\begin{equation}
U\equiv X_1^{-1}Y^{-1}Z^{-1}X_1 X_2 Y Z X_2^{-1}\label{eq:threevar}
\end{equation}
where $X_1 X_2$, $Y$, and $Z$ are nonempty, although one of $X_1$ or $X_2$ may 
be empty.  If \eqref{eq:twovar} holds, $U=[X_2 X_1,X_1^{-1}YX_2^{-1}]$ gives 
the first conclusion of Theorem \ref{th:BEF}, and 
$U^*=X_2^{-1}X_1^{-1}Y^{-1}X_1 X_2 Y=[X_1 X_2,Y]$ gives the second.  If 
\eqref{eq:threevar} obtains, $U=[X_2 Y X_1,X_1^{-1}Z X_2^{-1}]$ gives the 
first assertion of the theorem.  To see the second set of inequalities, we let 
$U^*=X_2^{-1}X_1^{-1}Y^{-1}Z^{-1}X_1 X_2 YZ$ and note that 
$2|X_1 X_2|+2|Y|+2|Z|=|U|$ implies that one of $|X_1 X_2|$, $|Y|$, $|Z|$ is 
less than or equal to $(1/6)|U|$.  If $|X_1 X_2|\le (1/6)|U|$, we write 
$U^*=[YX_1 X_2,X_2^{-1}X_1^{-1}Z]$; if $|Y|\le(1/6)|U|$, we write 
$U^*=[YX_1 X_2,YZ]$; if $|Z|\le(1/6)|U|$, we write $U^*=[Z^{-1}X_1 X_2,YZ]$.  
In each case, we see the truth of the second assertion of Theorem 
\ref{th:BEF}.
\end{pf*}

Theorem 1 of \cite{BEF} does not include the hypothesis that $U$ is cyclically 
reduced, and asserts that every solution to \eqref{eq:comm} is in the 
stabilizer class of a solution $\phi_0$ with $|x\phi_0|\le|U|-3$ and 
$|y\phi_0|\le|U|-3$.  We lose no generality in assuming that $U$ is cyclically 
reduced, for equations $[x,y]=C^{-1}U_0C$ and $[CxC^{-1},CyC^{-1}]=U_0$ have 
the same solutions.

The bounds given in Theorem \ref{th:BEF} for $|x\phi|+|y\phi|$ are not tight in 
the case that $|U|=4$, for then every solution is in the stabilizer class of 
a solution $\phi_0$ of \eqref{eq:comm} 
with $|x\phi_0|=|y\phi_0|=1$.  For $|U|>4$, though, we shall show that the 
bounds given in Theorem \ref{th:BEF} do not admit improvement.  

Let $n$ be a positive integer and let 
$\{a_1,\dots,a_n,b_1,\dots,b_n,c_1,\dots,c_n\}$ be a subset of $A$, our 
chosen set of free generators for $H$.  Consider first
\begin{displaymath}
U_1=b_n^{-1}\dots b_1^{-1}c_1^{-1}b_1\dots b_na_1\dots a_nc_1a_n^{-1}\dots 
a_1^{-1}.
\end{displaymath}
We find that, cyclically, $U_1$ is not a cancellation-free image of 
$x^{-1}y^{-1}z^{-1}xyz$, and that $U_1$ is a cancellation-free image of 
$x^{-1}y^{-1}xy$ in only one way, up to change of variables.  Thus every 
solution to $[x,y]=U_1$ is in the stabilizer class of $\phi_0$, where 
$x\phi_0=a_1\dots a_nb_1\dots b_n$ and 
$y\phi_0=a_1\dots a_nc_1a_n^{-1}\dots a_1^{-1}$.  Now elements of 
$\Stab([x,y])$ are automorphisms of $\langle x,y;\;\rangle$, so for any 
solution $\phi$ to \eqref{eq:comm}, $\{ x\phi,y\phi\}$ and 
$\{ x\phi_0,y\phi_0\}$ generate the same subgroup of $H$.  Since 
$\{ x\phi_0,y\phi_0\}$ is Nielsen reduced, it follows that 
for any solution $\phi$ to 
$[x,y]=U_1$, $|x\phi|+|y\phi|\ge|x\phi_0|+|y\phi_0|=|U_1|-1$.

Now let 
\begin{displaymath}
U_2=a_n^{-1}\dots a_1^{-1}b_n^{-1}\dots b_1^{-1}c_n^{-1}\dots c_1^{-1} 
        a_1\dots a_nb_1\dots b_nc_1\dots c_n.
\end{displaymath}
A typical cyclic permutation of $U_2$ is 
\begin{displaymath}
U_2^*=a_i^{-1}\dots a_1^{-1}b_n^{-1}\dots b_1^{-1}c_n^{-1}\dots c_1^{-1}
        a_1\dots a_n b_1\dots b_n c_1\dots c_n a_n^{-1}\dots a_{i+1}^{-1}
\end{displaymath}
with $0\le i\le n$.  As above, $[x,y]=U_2^*$ has only one stabilizer class of 
solutions, namely that of $\phi_0^*$ given by 
\begin{align*}
x\phi_0^*&=a_{i+1}\dots a_n b_1\dots b_n a_1\dots a_i ,\\
y\phi_0^*&=a_i^{-1}\dots a_1^{-1}c_1\dots c_n a_n^{-1}\dots a_{i+1}^{-1} .
\end{align*}
Again we find that $\{ x\phi_0^* , y\phi_0^* \}$ is Nielsen reduced, so for 
any solution $\phi^*$ to $[x,y]=U_2^*$, 
$|x\phi^*|+|y\phi^*|\ge|x\phi_0^*|+|y\phi_0^*|=(2/3)|U_2|$.

By similar arguments, we can see that for any cyclic permutation $U_1^*$ of 
$U_1$ and solution $\phi^*$ of $[x,y]=U_1^*$, either 
$|x\phi^*|\ge(1/2)|U_1|-1$ or $|y\phi^*|\ge(1/2)|U|-1$.

Our next result shows that an equation of the form $[x,y]=T$ with $T\in H$ may 
have any number of stabilizer classes of solutions.

\begin{theorem}\label{th:powers}
Suppose that $H=K\ast L$ with $K$ and $L$ free, that $U\in K$ and $V\in L$ are 
nontrivial and not proper powers, and that $m$ and $n$ are positive integers.  
The equation 
\begin{equation}\label{eq:powers}
[x,y]=[U^m,V^n]
\end{equation}
has $m+n-1$ distinct stabilizer classes of solutions represented by 
$x\phi_i=U^m$, $y\phi_i=U^{-i}V^n$ for $0\le i<m$ and by 
$x\psi_j=V^j U^m$, $y\psi_j=V^n$ for $1\le j<n$.
\end{theorem}

\begin{pf}
We note that if $\alpha$ is an automorphism of $H$, solutions $\mu$ and $\nu$ 
of \eqref{eq:powers} are in the same stabilizer class if and only if 
$\mu\alpha$ and $\nu\alpha$ are in the same stabilizer class of solutions to 
$[x,y]=[(U\alpha)^m,(V\alpha)^n]$.  Applying an automorphism of $H$ induced by 
inner automorphisms of $K$ and $L$, then, we may assume that $U$ and $V$ are 
cyclically reduced.

We next need to find the ways in which $[U^m,V^n]$, as a cyclic word, can be 
factored as $X^{-1}Y^{-1}XY$ or $X^{-1}Y^{-1}Z^{-1}XYZ$.  In any such 
factorization, $X^{\pm 1}$, $Y^{\pm 1}$, and $Z^{\pm 1}$ must be subwords of 
$U^{\pm m}$ and $V^{\pm n}$.  If, for example, $X^{-1}$ overlapped both 
factors of the product $U^{-m}V^{-n}$, then $X$ would overlap both factors of 
either $V^{-n}U^m$ or $V^nU^{-m}$, but this would violate the cyclic reduction 
of either $V$ or $U$.  If $[U^m,V^n]=X^{-1}Y^{-1}XY$ as cyclic words, then, it 
must be that up to change of notation $X=U^m$ and $Y=V^n$; this gives solution 
$\phi_0$.  If $[U^m,V^n]=X^{-1}Y^{-1}Z^{-1}XYZ$ as cyclic words, one of $X$, 
$Y$, $Z$ must be $U^{\pm m}$ or $V^{\pm n}$.  If, say, $Z=V^n$, then 
$XY=U^m$ and $X^{-1}Y^{-1}=U^{-m}$, $X$ and $Y$ commute, and so $X$, $Y$, and 
$U$ are powers of a common element.  Since $U$ is not a proper power, 
$X=U^i$ and $Y=U^{m-i}$ for some $i$, $1\le i<m$.  This gives us the solution 
$x\mapsto Z^{-1}X=V^{-n}U^i$, $y\mapsto YZ=U^{m-i}V^m$, which is in the 
stabilizer class of $\phi_i$.  Likewise if, say, $X=U^m$, then $Y=V^j$ and 
$Z=V^{n-j}$ for some $j$, $1\le j<n$, which gives rise to the solution 
$x\mapsto Z^{-1}X=V^{-n+j}U^m$, $y\mapsto YZ=V^n$, in the stabilizer class of 
$\psi_j$.

It remains to show that 
$\phi_0,\phi_1,\dots,\phi_{m-1},\psi_1,\dots,\psi_{n-1}$ are 
in different stabilizer classes.  We do this by proving that the subgroups 
$S_i$, $T_j$ of $H$ generated by $\{x\phi_i,y\phi_i\}$ for $0\le i<m$ and by
$\{x\psi_j,y\psi_j\}$ for $1\le j<n$ are all different.

We first show that if $0\le i<j<m$, then $S_i\ne S_j$.  We do this by cases.  
First, if $|V^n|\le|U^i|$ and $|V^n|\le|U^{m-i}|$, $\{U^{-i}V^n,U^{m-i}V^n\}$ 
is a Nielsen reduced generating set for $S_i$ and we see that 
$U^{m-j}V^n\in S_j - S_i$.  Next, if $i\le m/2$ and $|V^n|>|U^i|$ or 
$|V^n|>|U^{m-i}|$, then $|V^n|>|U^i|$ and $\{U^{-i}V^n,U^m\}$ is a Nielsen 
reduced generating set for $S_i$.  We again note that 
$U^{m-j}V^n\in S_j -S_i$.  
Finally, if $i>m/2$ and either $|V^n|>|U^i|$ or $|V^n|>|U^{m-i}|$, then 
$|V^n|>|U^{m-i}|$ and $\{U^{m-i}V^n,U^m\}$ is a Nielsen reduced generating set 
for $S_i$, and once more $U^{m-j}V^n\in S_j - S_i$.

In a similar way, we show that for $1\le i<j<n$, $T_i\ne T_j$.  Finally, we 
note that $U^m$ is a member of each of the $S_i$ but none of the $T_j$, which 
distinguishes the $S_i$ from the $T_j$.
\end{pf}

\begin{theorem}\label{th:squares}
If $a_1,\dots a_g$ are distinct free generators of a free group $H$ and if 
$n_1,\dots n_g$ are nonzero integers, every solution to 
\begin{equation}\label{eq:squares}
x_1^2\dots x_g^2=a_1^{2n_1}\dots a_g^{2n_g}
\end{equation}
is in the stabilizer class of the solution $\phi_0$ given by 
$x_i\phi_0=a_i^{n_i}$ for $1\le i\le g$.
\end{theorem}

\begin{pf}
Let $U=a_1^{2n_1}\dots a_g^{2n_g}$.  It is clear that $U\not\in H'$, that 
$U\in 2H$, and that $\genus^-(U)\le g$.  Suppose that $\genus^-(U)=k$ and let 
$\phi$ be a solution to $x_1^2\dots x_k^2=U$.  By Theorem \ref{th:main}, there 
is an irredundant nonorientable quadratic $W\in F=\langle x_1,\dots;\;\rangle$ 
with $\genus^-(W)=k$ and a cancellation-free solution $\psi$ to $W=U$ such 
that if $\gamma_W$ is an automorphism of $F$ with 
$W\gamma_W=x_1^2\dots x_k^2$, then $\phi=\sigma\gamma_W^{-1}\psi$ for some 
$\sigma\in\Stab(x_1^2\dots x_k^2)$.  

There is no variable $y$ such that both $y$ and $y^{-1}$ occur in $W$, for 
there is no $Y\in H$, $Y\ne 1$, such that both $Y$ and $Y^{-1}$ are 
(cyclically) subwords of $U$.  Further, the form of $U$ shows that for each 
letter $y$ in $W$, $y\psi$ is contained in a single syllable $a_i^{2n_i}$
of $U$.  Thus $W\equiv W_1\dots W_g$ with $W_1,\dots W_g$ words on disjoint 
sets of variables and with $W_i\psi=a_i^{2n_i}$ for $1\le i\le g$.  It follows 
that each $W_i$ is nonorientable quadratic with $\genus^-(W_i)\ge 1$, so 
$k=\genus^-(W)=\sum_{i=1}^g \genus^-(W_i)\ge g$.  Thus $\genus^-(W)=g$ and 
$\genus^-(W_i)=1$ for $1\le i\le g$, so each $W_i$ is the square of a variable 
or its inverse.  We may take $\gamma_W$ to be an automorphism of $F$ that 
permutes $\{ x_1,x_1^{-1},\dots\}$, and so $\gamma_W^{-1}\psi=\phi_0$ and 
$\phi$ is in the stabilizer class of $\phi_0$.
\end{pf}

\section{Commutators as Products of Squares}\label{se:comms}

It is well known that in a free group, indeed in any group, every commutator 
is a product of three squares:
\begin{displaymath}
[U,V]=(U^{-1})^2 (UV^{-1})^2 V^2.
\end{displaymath}
One cannot in general get by with fewer than three squares; for example, 
Lyndon and Morris Newman \cite{LynNew73} have shown that in the free group on 
$a$ and $b$, $[a,b]$ is not a product of two squares.  We can see that this is 
true, for $[a,b]$ is not a cancellation-free image of any of the nonorientable 
Wicks forms of genus two.

In a free group, a nontrivial commutator is never a square; this was first 
noted by M.~Sch\"{u}tzenberger \cite{Schutz59}, and follows from examination 
of the orientable Wicks forms of genus one.  Thus one is left with the 
question of which commutators are products of two squares in a free group.  
One obvious possibility is that $[U,V]=[S^2,T]$ for some $S$ and $T$, in which 
case $[U,V]=(S^{-1})^2 (T^{-1}ST)^2$.  This includes the situation that 
$[U,V]=[P,Q^2]$ for some $P$ and $Q$, for $[P,Q^2]=[P^{-1}Q^2P,P^{-1}]$. One 
might ask if this is the only way in which a commutator can be a product of 
two squares in a free group.  We shall give an example to show that this is not 
the case.  This complements an example given by J.~Comerford and Y.~Lee 
\cite{ComLee} to show that if in a free group a product of two commutators is 
a square, it need not be the square of a commutator.

Our example fits into a sequence of results about solutions to equations of 
the form
\begin{equation}\label{eq:homog}
x_1^2\dots x_g^2=1
\end{equation}
in a free group.  It is easy to see that if $g=1$, \eqref{eq:homog} implies 
that $x_1=1$, and that if $g=2$, \eqref{eq:homog} implies that $x_2=x_1^{-1}$.  
For $g=3$, Lyndon proved \cite{Lyn59} that \eqref{eq:homog} implies that 
$x_1$, $x_2$, and $x_3$ are powers of a common element.  Thus for $g\le3$, 
\eqref{eq:homog} has only \lq\lq obvious\rq\rq\ solutions.  When $g=5$, 
\eqref{eq:homog} can be rewritten using an automorphism of 
$\langle x_1,\dots,x_5;\;\rangle$ as $x_1^2=[x_2,x_3][x_4,x_5]$; J.~Comerford 
and Lee showed \cite{ComLee} that this has \lq\lq nonobvious\rq\rq\ solutions.  
The case we consider here is $g=4$, in which \eqref{eq:homog} can be rewritten 
as $x_1^2x_2^2=[x_3,x_4]$.  Again, we show that this has 
\lq\lq nonobvious\rq\rq\ solutions. 

We now give our example.  
Let $U=b^{-1}a^{-1}b^2ab^{-1}$ and $V=a$ in $H=\langle a,b;\;\rangle$.  One 
can check that 
\begin{displaymath}
[U,V]=(ba^{-1}b^{-1}a^{-1}b^{-1}abab^{-1})^2 
(ba^{-1}b^{-1}a^{-1}b^2ab^{-1}a)^2,
\end{displaymath}
but we shall show that $[x^2,y]=[U,V]$ has no solutions in $H$.  This follows 
from the following two results.

\begin{lemma}\label{le:comm}
Every solution to $[x,y]=[U,V]$ is in the stabilizer class of $x\phi_0=U$, 
$y\phi_0=V$.
\end{lemma}

\begin{lemma}\label{le:subgp}
If $K$ is the subgroup of $H$ generated by $U$ and $V$, if $W\in K$, and if 
$x^2=W$ has a solution in $H$, then $x^2=W$ has a solution in $K$.
\end{lemma}

To see that these lemmas imply that $[x^2,y]=[U,V]$ has no solutions, suppose 
that $\psi$ is a solution to $[x^2,y]=[U,V]$.  Let $x\psi=P$ and $y\psi=Q$.  
By Lemma \ref{le:comm}, $P^2=x\sigma\phi_0$ and $Q=y\sigma\phi_0$ for some 
$\sigma\in\Stab([x,y])$.  Thus $P^2$ and $Q$ are elements of $K$ and, by Lemma 
\ref{le:subgp}, $P$ is an element of $K$.  Now $U$ and $V$ are free generators 
for $K$, so $\beta : K\to\langle x,y;\;\rangle$ defined by $U\beta=x$ and 
$V\beta=y$ is an isomorphism.  It follows that $\beta\sigma\phi_0$ is an 
automorphism of $K$ that sends $U$ to $P^2$, so $U$ is the square of an 
element of $K$, which is impossible.

\begin{pf*}{Proof of Lemma \ref{le:comm}}
We must show that as a cyclic word,
\begin{displaymath}
[U,V]=ba^{-1}b^{-2}ab a^{-1} b^{-1}a^{-1}b^2ab^{-1} a
\end{displaymath}
does not factor as $X^{-1}Y^{-1}Z^{-1}XYZ$ and factors as 
$X^{-1}Y^{-1}XY$ in essentially only one way.  

Suppose first that there is a factorization of the cyclic word $[U,V]$ in 
which the copy of $ab$ consisting of the last and first letters of the 
ordinary word $[U,V]$ is contained in one of 
$X^{\pm 1}$, $Y^{\pm 1}$, $Z^{\pm 1}$, say in $X$.  Then 
$X\equiv X_1 abX_2$ with $X_1$ and $X_2$ possibly empty, and since 
$b^{-1}a^{-1}$ occurs in $[U,V]$ in only one position, 
\begin{displaymath}
a^{-1}b^{-2}aba^{-1}=X_2YX_2^{-1}\text{ or }X_2YZX_2^{-1}
\end{displaymath}
and
\begin{displaymath}
b^2ab^{-1}=X_1^{-1}Y^{-1}X_1\text{ or }X_1^{-1}Y^{-1}Z^{-1}X_1.
\end{displaymath}
But these imply that $a^{-1}b^{-2}aba^{-1}$ is conjugate to 
$(b^2ab^{-1})^{-1}$, which is not the case.

It must be, then, that up to change of notation,
\begin{displaymath}
ba^{-1}b^{-2}aba^{-1}b^{-1}a^{-1}b^2ab^{-1}a=X^{-1}Y^{-1}XY\text{ or }
X^{-1}Y^{-1}Z^{-1}XYZ.
\end{displaymath}
By a length comparison, 
\begin{displaymath}
bab^{-2}aba^{-1}=X^{-1}Y^{-1}\text{ or }X^{-1}Y^{-1}Z^{-1} 
\end{displaymath}
and
\begin{displaymath}
b^{-1}a^{-1}b^2ab^{-1}a=XY\text{ or }XYZ.
\end{displaymath}
Thus $X$ begins and ends with $b^{-1}$.  The possibilities are 
$X=b^{-1}a^{-1}b^2ab^{-1}=U$ and $Y=a=V$ or $X=b^{-1}$.  But $X=b^{-1}$ 
implies that $a^{-1}b^2ab^{-1}a$ is conjugate to $(ab^{-2}aba^{-1})^{-1}$, 
which isn't so.
\end{pf*}

\begin{pf*}{Proof of Lemma \ref{le:subgp}}
Let us suppose that $W\ne 1$.  Since $\{U,V\}$ is a free generating set for 
$K$, there is a unique expression $W=Z(U,V)$ with $Z$ a freely reduced word in 
$U$ and $V$.  In fact, the map $x\mapsto U$, $y\mapsto V$ is a 
cancellation-free solution to $Z(x,y)=W$.  We may assume without loss of 
generality that $Z(U,V)$ is a cyclically reduced word on $U$ and $V$, and 
hence that $W$ is a cyclically reduced word on $a$ and $b$.

Now suppose that $W=T^2$ for some (cyclically reduced) $T\in H$.  If 
$T\not\in K$, then $T\equiv T_0U_1$ with $T_0\in K$ and $U=U_1U_2$ a 
nontrivial factorization of $U=b^{-1}a^{-1}b^2ab^{-1}$.  But this implies that 
some nontrivial initial subword $U_1$ of $U$ is an element of $K$, which is 
plainly not the case.  Thus $T\in K$ and $W$ is a square in $K$.
\end{pf*}

We close with a brief description of how we found this example.  Using 
generators for the nonorientable mapping class groups given by J.~Birman and 
D.~Chillingworth \cite{BirChil}, J.~Comerford and Lee \cite{ComLee} provided 
generating sets for the stabilizers of nonorientable quadratic words.  They 
gave a set of five generators, $\Bar{a}_1$, $\Bar{b}_1$, $\Bar{c}_1$, 
$\Bar{b}_2$, and $y_3$, for the stabilizer of 
$W=s_1t_1s_1^{-1}t_1^{-1}s_2t_2s_2^{-1}t_2$.  Now if $\alpha$ is the 
automorphism of $\langle s_1,t_1,s_2,t_2;\;\rangle$ that sends $t_2$ to 
$s_2t_2$ and fixes $s_1$, $t_1$, and $s_2$, we see that 
$W\alpha=s_1t_1s_1^{-1}t_1^{-1}s_2^2t_2^2$, and that $\sigma\in\Stab(W)$ if 
and only if $\alpha^{-1}\sigma\alpha\in\Stab(W\alpha)$.  The idea, then, is to 
start with the \lq\lq obvious\rq\rq\ solution $\phi_0$ to $W\alpha=1$ given by 
$s_1\phi_0=1$, $t_1\phi_0=a$, $s_2\phi_0=b^{-1}$, $t_2\phi_0=b$, and preceed 
it with elements of $\Stab(W\alpha)$ to produce \lq\lq nonobvious\rq\rq\ 
solutions to $W\alpha=1$.  The example presented here is based on the solution 
$\alpha^{-1}y_3^2\alpha\phi_0$ to $W\alpha=1$.

\newpage

\ifx\undefined\bysame
\newcommand{\bysame}{\leavevmode\hbox to3em{\hrulefill}\,}
\fi

\end{document}